\documentclass[12pt]{article}

\usepackage{amsmath}
\usepackage{amssymb}
\usepackage{amsthm}
\usepackage{a4wide}
\usepackage[english]{babel}

\def\Dem{\hspace{-\parindent}{\em Proof }}

\theoremstyle{plain}
\newtheorem{Th}{Theorem}[section]
\newtheorem{Prop}[Th]{Proposition}

\newtheorem{Lemme}[Th]{Lemma}

\theoremstyle{definition}

\newtheorem{Rem}[Th]{Remark}

\newcommand{\N}{\mathbb{N}}
\newcommand{\Z}{\mathbb{Z}}
\newcommand{\Q}{\mathbb{Q}}
\newcommand{\R}{\mathbb{R}}
\newcommand{\C}{\mathbb{C}}

\newcommand{\drondsoul}{\underline{\drond}}
\newcommand{\opnv}{{\rm Op}}
\newcommand{\opW}{\opnv_W}
\newcommand{\opWTH}{\opnv_{W \inter T_0(H)}}
\newcommand{\opWTHpr}{\opnv_{W \inter T_0(H')}}
\newcommand{\opWTHsec}{\opnv_{W \inter T_0(H'')}}
\newcommand{\opdrondT}{\opnv _{\drondsoul, \Tsoul}}
\newcommand{\Ndt}{\N ^d _{\Tsoul}}
\newcommand{\Tsoul}{\underline{T}}
\newcommand{\Span}{{\rm Span}}
\newcommand{\muet}{\mu^*}
\newcommand{\inclus}{\subset}
\newcommand{\inter}{\cap}
\newcommand{\drond}{\partial}
\newcommand{\oo}{{\cal O}}
\renewcommand{\P}{{\mathbb P}}
\newcommand{\rk}{{\rm rk }}
\newcommand{\Card}{{\rm Card }}
\newcommand{\eps}{\varepsilon}

\newcommand{\und}{\{1, \ldots, d\}}

\newcommand{\unr}{\{1, \ldots, r\}}

\newcommand{\zeror}{\{0,\ldots,r\}}

\newcommand{\zerormu}{\{0,\ldots,r-1\}}
\newcommand{\unrmu}{\{1,\ldots,r-1\}}
\newcommand{\unlmu}{\{1, \ldots, l-1 \}}   
\newcommand{\unn}{\{1, \ldots, n\}}

\newcommand{\cale}{{\mathcal E}}
\newcommand{\calE}{{\mathcal E}}
\newcommand{\calS}{{\mathcal S}}

\newcommand{\Sgothhprhsec}{{\mathcal N}_{H',H''}(\Ssoul)}
\newcommand{\SgothKzero}{{\mathcal N}_{K,\{0\}}(\Ssoul)}

\newcommand{\chaine}{(H_i)_{0 \leq i \leq r}}
\newcommand{\sss}{{\mathfrak S}}

\newcommand{\base}{{\mathcal B}_{G, \GS, T, D}}
\newcommand{\basesanssoul}{{\mathcal B}_{G, \Gamma(S), T, D}}
\newcommand{\baseom}{{\mathcal B}_{G, \Omega, T, D}}
\newcommand{\basek}{{\mathcal B}_{G, \GS, k \alpha , k}}

\newcommand{\Gbar}{{\overline{G}}}

\newcommand{\cdegD}{{\cal R}(G) _{D}} 
\newcommand{\cdeg}{{\cal R}(G)} 

\newcommand{\cdehi}{{\cal R}(H_i)} 
\newcommand{\cdequo}{{\cal R}(H_{i+1}/H_i)} 
\newcommand{\cdehipu}{{\cal R}(H_{i+1})} 
\newcommand{\cdegsh}{{\cal R}(G/H_{i+1})}

\newcommand{\xzd}{X_0^D}

\newcommand{\Sul}{S_1,\ldots,S_l}
\newcommand{\gamul}{\gamma_1,\ldots,\gamma_l}
\newcommand{\Sulgeq}{S_1 \geq \ldots \geq S_l}
\newcommand{\Ssoul}{\underline{S}}
\newcommand{\GS}{\Gamma(\Ssoul)}
\newcommand{\glmu}{\Gamma_{l-1}}
\newcommand{\gti}{\widetilde\gamma}

\newcommand{\cas}{connected algebraic subgroup }
\newcommand{\cass}{connected algebraic subgroups }
\newcommand{\fis}{\varphi_{\Ssoul}}
\newcommand{\fist}{\varphi_{\Ssoul, \Tsoul}}
\newcommand{\fishi}{\fis(H_i)}
\newcommand{\fishipu}{\fis(H_{i+1})}
\newcommand{\fishimu}{\fis(H_{i-1})}
\newcommand{\dimhi}{\dim H_i }
\newcommand{\dimhipu}{\dim H_{i+1} }
\newcommand{\dimhimu}{\dim H_{i-1} }

\renewcommand{\hbar}{\overline{h}}
\newcommand{\Hbar}{\overline{H}}
\newcommand{\betabar}{\overline{\beta}}
\newcommand{\Omegabar}{\overline{\Omega}}

\newcommand{\cintrou}{c_1}
\newcommand{\cintrod}{c_2}
\newcommand{\cthintro}{c_3}
\newcommand{\cthuu}{c_4}
\newcommand{\cthud}{c_5}
\newcommand{\cconj}{c_6}
\newcommand{\cpphz}{c_7}
\newcommand{\clemfn}{c_8}
\newcommand{\cinterp}{c_9}
\newcommand{\caba}{c_{10}}
\newcommand{\cabb}{c_{11}}
\newcommand{\cabc}{c_{12}}
\newcommand{\ctss}{c_{13}}
\newcommand{\cnvu}{c_{14}}
\newcommand{\cnvt}{c_{15}}
\newcommand{\cabbpr}{c_{16}}
\newcommand{\ctst}{c_{17}}
\newcommand{\cnvd}{c_{18}}
\newcommand{\cabapr}{c_{19}}

\title{Connecting Interpolation and Multiplicity Estimates in Commutative Algebraic Groups}
\author{S. Fischler and M. Nakamaye} 
\date{\today}

\begin{document}

\maketitle

\begin{abstract}
Let $G$ be a  commutative algebraic group  embedded in projective space and $\Gamma$ a finitely generated subgroup of $G$. From these data we construct a chain of algebraic subgroups of $G$ which is intimately related to obstructions to multiplicity or interpolation estimates. Let  $\gamma_1,\ldots,\gamma_l$ denote a family of generators of $\Gamma$ and, for any $S>1$, let  $\Gamma(S)$  be  the set of  elements $n_1\gamma_1+\ldots+n_l\gamma_l$ with integers $n_j$ such that $|n_j| < S$. Then this chain of subgroups controls, for large values of $S$, the distribution of $\Gamma(S)$ with respect to algebraic subgroups of $G$. As an application we essentially determine (up to multiplicative constants) the locus of common zeros of all $P \in  H^0(\overline G ,{\cal O}(D))$ which vanish to at least  some given order at all points of $\Gamma(S)$. When $D$ is very small this result reduces to a multiplicity estimate; when $D$ is very large it is a kind of interpolation estimate.
\end{abstract}

\bigskip

\noindent{\bf Math. Subj. Classification (2010):} 
14L10 (Group varieties); 
11J95   (Results involving abelian varieties);
14L40   (Other algebraic groups (geometric aspects));
 14C20 (Divisors, linear systems, invertible sheaves).

\section{Introduction} \label{sec0}

Let $G$ 
be a positive dimensional connected commutative algebraic group,  embedded in $\P^N$ through the choice of a very ample divisor on a 
compactification of $G$.  In most transcendence proofs involving $G$, an important role is played by the evaluation map
\begin{equation} \label{eqeval}
 H^0(\Gbar ,\oo(D)) \to  H^0\left(\Gbar , \oo(D) \otimes \oplus_{\omega \in \Gamma(S)}\oo_\Gbar/m_\omega^{T}  \right) ; 
 \end{equation}
here $\Gbar$ is the Zariski closure of $G$ in $\P^N$, $m_\omega \subset \oo_\Gbar$  is the maximal ideal sheaf corresponding to the point $\omega$,  and for a positive real number $S$,  
$\Gamma(S)$ is the set of all elements $n_1\gamma_1+\ldots+n_l\gamma_l$ with integers $n_j$ such that $|n_j| < S$. In this setting $\gamul$ are fixed elements of $G$  and $S$ is often chosen to be very large.
We let $\Gamma$ denote the  $\Z$-module generated by $\gamul$.  The set  $\Gamma(S)$ depends 
on $\gamul\in G$ in addition to  $\Gamma$ and $S$, making the notation  $\Gamma(S)$
 rather unpleasant but it   is the usual one in this setting. The integers $D$, $S$, $T$ are parameters which typically take very large values in transcendence proofs, except when no multiplicities are involved, that is when  $T=1$. 

A crucial step in most transcendence proofs is the {\em multiplicity estimate}, called a {\em zero estimate} when $T=1$. The simplest one in this setting is perhaps the following (see \cite{PPhzeros} or \cite{Wustholz1989}):   if $D < \cintrou TS^\mu$ then \eqref{eqeval} is injective so that $P=0$ as soon as $P \in  H^0(\Gbar ,\oo(D)) $ vanishes to order at least $T$ at each point $\omega\in\Gamma(S)$. Here $\cintrou$ is a positive constant depending on $G$, its embedding in $\P^N$, and $\gamul$. The real exponent $\mu\geq 0$ is defined by
$$\mu = \mu(\Gamma,G) = \min_{H \subsetneq G} \frac{\rk(\Gamma)-\rk(\Gamma\cap H)}{\dim G - \dim H}$$
where $H$ ranges through the set of all proper \cass of $G$.

Instead of a multiplicity estimate and the construction of an auxiliary function, it is possible to use an {\em interpolation estimate} and an auxiliary functional (see \cite{MiWDelphes}, \cite{MiWFB}, \cite{Grundlehren}, \cite{MiWdep}). Such a result was proved by Masser \cite{Masser} when no multiplicities are involved, that is when  $T=1$, and generalized by the first author \cite{SFinterp}. It reads as follows: if $D > \cintrod TS^{\muet}$ then \eqref{eqeval} is surjective, where 
$\cintrod$ is a positive constant depending on $G$, its embedding in $\P^N$, and $\gamul$. The real exponent $\muet\geq 0$ is defined by
$$\muet = \muet(\Gamma,G) = \max_{H \neq \{0\}} \frac{ \rk(\Gamma\cap H)}{  \dim H}$$
where $H$ ranges through the set of all non-zero \cass of $G$.

The exponents $ \mu(\Gamma,G)$ and $ \muet(\Gamma,G)$ measure the distribution of $\Gamma$ (and that of $\Gamma(S)$, if  $S$ is sufficiently large) with respect to algebraic subgroups of $G$. The former appears in early zero estimates \cite{MWun} and already in  \cite{MiWAsterisque} (\S 1.3). It is related to the density coefficient of $\Gamma$ if $G = \mathbb{G}_a^n$ and 
$\Gamma \subset (\overline{\Q} \cap \R)^n$ (see \S 1.3.d of   \cite{MiWAsterisque}),    and to Schwarz lemmas (see Chapter 7 of  \cite{MiWAsterisque} and \cite{RoySchwarz}). The exponent $ \muet(\Gamma,G)$ is a dual version introduced in  \cite{Masser}.  These exponents satisfy the inequalities
$$\mu(\Gamma,G) \leq \frac{\rk \, \Gamma}{\dim G} \leq  \muet(\Gamma,G)$$
by definition.  A finitely generated $\Z$-module $\Gamma\subset G$ is said to be {\em well distributed} in $G$ if $\mu(\Gamma,G) = \frac{\rk \, \Gamma}{\dim G}$ or,  equivalently, if $\muet(\Gamma,G) = \frac{\rk \, \Gamma}{\dim G}$ (see \cite{Masser}).

Elaborating upon ideas of \cite{NakamayeLang}, we construct in \S \ref{sec3} a chain of algebraic subgroups $\{0\} = H_0 \subsetneq H_1  \subsetneq \ldots   \subsetneq H_r = G$, with $r\geq 1$, associated 
to $\Gamma$ and $G$. These subgroups satisfy
$$\mu(\Gamma \cap H_j \bmod H_i , H_j/H_i) = \frac{\rk(\Gamma\cap H_j) - \rk(\Gamma\cap H_{j-1})}{\dim H_j - \dim H_{j-1}}$$
and
$$\muet(\Gamma \cap H_j \bmod H_i , H_j/H_i) = \frac{\rk\Big(\frac{ \Gamma\cap H_{i+1}}{ \Gamma\cap H_i}\Big)}{\dim ( H_{i+1}/H_i)}$$
for any $i$, $j$ such that $0\leq i < j \leq r$. Here and throughout this text, we let $\Omega \bmod H  = \frac{\Omega+H}{H}$ for any subset $\Omega$ of $G$ and any algebraic subroup $H \subset G$.

The following properties hold:
\begin{itemize}
\item $ \mu(\Gamma,G)  =  \frac{\rk(\Gamma  ) - \rk(\Gamma\cap H_{r-1})}{\dim  G  - \dim H_{r-1}}$.
\item $ \muet(\Gamma,G)  =  \frac{ \rk(\Gamma\cap H_{ 1})}{   \dim H_{1}}$.
\item $\Gamma$  is  well distributed in $G$ if and only if $r=1$ so that the chain is simply $\{0\} = H_0 \subsetneq H_1 = G$.
\item For any $i\in\zerormu$, $\Gamma\cap H_{i+1}\bmod  H_i$ is well distributed in $H_{i+1}/H_i$.
\end{itemize}
Moreover, if $H$ is a non-zero \cas of $G$ such that $ \muet(\Gamma,G) =   \frac{ \rk(\Gamma\cap H)}{  \dim H}$, then $H \subset H_1$ (see \cite{NakamayeLang}, \S 1.3). In the same way,  if $H$ is a proper  \cas of $G$ such that $ \mu (\Gamma,G) =  \frac{\rk(\Gamma)-\rk(\Gamma\cap H)}{\dim G - \dim H}$, then $H_{r-1} \subset H $.

We think this chain of subgroups can be useful in many problems where the distribution of $\Gamma$  with respect to algebraic subgroups of $G$ is involved, for instance in studying
 the points of $\Gamma(S)$ in the spirit of  \S 1.3.d of   \cite{MiWAsterisque} for the
case  $G = \mathbb{G}_a^n$ and $\Gamma \subset (\overline{\Q} \cap \R)^n$. It may also provide
 a geometric interpretation
 closely related to the Seshadri exceptional subvarieties studied in \cite{FNinterp}.
  We use these subgroups here to study the locus $\basesanssoul$ of common zeros of all $P$ in the kernel of \eqref{eqeval}, that is the set of $x\in G$ such that $P(x) = 0$ for any $P \in  H^0(\Gbar ,\oo(D)) $ which vanishes to order at least $T$ at each point of $\Gamma(S)$. 

\bigskip

To state our result, we let
$$\mu_i  = \mu(\Gamma \cap H_{i+1} \bmod H_i , H_{i+1}/H_i) = \frac{\rk(\Gamma\cap H_{i+1}) - \rk(\Gamma\cap H_{i})}{\dim H_{i+1} - \dim H_i}$$
for any $i\in\zerormu$. As stated previously, we have also $\mu_i  = \muet(\Gamma \cap H_{i+1} \bmod H_i , H_{i+1}/H_i) $ since $\Gamma\cap H_{i+1}\bmod  H_i$ is well distributed in $H_{i+1}/H_i$. Moreover $\mu_0 = \muet(\Gamma,G)$ and  $\mu_{r-1} = \mu(\Gamma,G)$; we shall prove that 
\begin{equation} \label{eqmudecr}
\mu_0 > \mu_1 > \ldots > \mu_{r-1}. 
\end{equation}
For convenience we write $\mu_{-1} = +\infty$ and $\mu_r = -\infty$.
In loose terms, the series of inequalities \eqref{eqmudecr} can be understood as follows. 
The algebraic subgroup $H_1$ contains the largest possible proportion of $\Gamma$ (with respect to its dimension), so that  the  proportion of $\Gamma \bmod H_1$ contained in $H_2/H_1$ has to be smaller 
with respect to $\dim(H_2/H_1)$; otherwise $H_2$ would contradict the maximality of $H_1$. 
This argument is made precise in Proposition \ref{propchaineiii}  below and the associated remarks.

Our main result reads as follows.

\begin{Th} \label{thintro} For any $\eps> 0$ with $0  < \eps < 1$  there exists a positive constant $\cthintro$, depending only on the embedding
of $G$ in $\P^N$, $\eps$, and  $\gamul$, with the following property.
For any positive integers $D$ and $T$, if $S$ is a sufficiently large positive integer (in terms of $G \hookrightarrow \P^N$, $\eps$, $\gamul$) and
$$\cthintro S^{\mu_i} T \leq D \leq \cthintro^{-1} S^{\mu_{i-1}} T $$
for some $i\in\zeror$, then we have 
$$\Gamma((1-\eps)S)+H_i\subset \basesanssoul \subset \Gamma(S) + H_i.$$
\end{Th}

With $i=r$ this is the above-mentioned multiplicity estimate because $\basesanssoul = G$.  With $i=0$ it follows from an interpolation estimate
since such an  estimate gives sections which separate jets at the points of $\Gamma(S) \cup \{ x\}$ for any $x\not\in\Gamma(S)$.  This result establishes a bridge between multiplicity and interpolation  estimates. It is a partial answer to a question asked by Michel Waldschmidt to the first author: what can be said about the evaluation map \eqref{eqeval} if $D$ is too large to apply a multiplicity estimate but too small to apply an  interpolation estimate? Of course this question remains largely open: for instance no non-trivial lower bound on the rank of this linear map is known for these values of $D$.  However we hope that Theorem \ref{thintro} can be useful to produce new transcendence proofs. 

If $X$ is a smooth projective variety, $\eta \in X$ a very general point, and
$L$ an ample line bundle on $X$ then the analogue of \eqref{eqeval} has
been studied closely  (see \cite{nv1} and \cite{nv2}):
$$
H^0\left(X,L^{\otimes D}\right) \longrightarrow H^0\left(X,L^{\otimes D} \otimes {\cal O}_X/m_\eta^T\right).
$$
The main idea is that once $\frac{T}{D}$ excedes the Seshadri constant of $L$ at $\eta$,
then the map ceases to be surjective.  This failure is estimated in  \cite{nv1} and \cite{nv2},  and it is
this extra information which allows a quantitative improvement for the lower bound of the
Seshadri constant of $L$ at $\eta$.  These techniques have been formalized in a broader
setting in   \cite{nv3}.

\bigskip

 Another motivation for Theorem \ref{thintro} is its relation to a conjecture of the second author (see \S \ref{subsecconj}).   In
Conjecture 1.1.9 of \cite{NakamayeLang}   a sequence of subgroups analagous to our $\chaine$ is  alluded to 
  and it is conjectured that these
  subgroups appear as the base locus of a linear series as in Theorem \ref{thintro}.   Because the methods employed
  in that paper are restricted to working on a compactification of $G$, with no auxiliary constructions such as
  projections to quotient groups, it was not
  possible to bound from above the size of the base loci in question as is done here.

When $D$ lies between $ \cthintro^{-1} S^{\mu_{i }} T $  and  $ \cthintro S^{\mu_{i }} T $,  for some $i \in \zerormu$, Theorem \ref{thintro} applied with these bounds  yields
$$\Gamma((1-\eps)S)+H_i\subset \basesanssoul \subset \Gamma(S) + H_{i+1} $$
since $\basesanssoul$ is a non-increasing function of $D$ when the subset $\Gamma(S)$ and the order of vanishing $T$ are held constant.   It would be interesting to have more information on $\basesanssoul$ for these critical values of $D$, but new ideas are needed. Indeed the proof of  Theorem \ref{thintro} is based on applying the special cases $i=0$ and $i=r$ to sub-quotients of $G$ obtained from the chain of subgroups $\chaine$. This strategy, remniscent of that used by Masser \cite{Masser} to prove his interpolation estimate, is responsible for the constant $\cthintro$.

\bigskip

In this paper we shall prove   Theorem \ref{thintro}  in a more general form: for any $\Sul \in\R$ we consider the set $\GS$ of all points $n_1\gamma_1+\ldots+n_l\gamma_l$ with integers $n_j$ such that $|n_j| < S_j$. Here $\Ssoul$ denotes the tuple $(\Sul)$, and we let  $\lambda \Ssoul = (\lambda S_1,\ldots,\lambda S_l)$ for any $\lambda >0$. Up to a permutation of $\gamul$, we may assume that $\Sulgeq$. This assumption will be useful to define the subgroups $H_i$ which depend in this case on $\Sul$ and $\gamul$ (whereas they depend only on $\Gamma$ and $G$ if $S_1=\ldots=S_l$).  The distribution of $\GS$   with respect to algebraic subgroups of $G$ is no longer measured simply by exponents like $\mu$, $\muet$ and the $\mu_i$ (see for instance \S 3 of \cite{SFinterp}).

  For any subset $\Omega$ of $G$, we let $\baseom$ denote the set of $x\in G$ such that $P(x) = 0$ for any $P \in  \cdegD $ which vanishes to order at least $T$ at each point of $\Omega$; here and throughout this text, we let $\cdegD =   H^0(\Gbar ,\oo(D)) $ as soon as $G$ is a commutative algebraic group embedded in a projective space, and we call {\em homogeneous polynomial of degree $D$}  any element of $\cdegD$.
The base field  is $\C$, though any algebraically closed field of
characteristic zero could be considered, for instance its 
$p$-adic analog $\C_p$; see also \cite{Masser}, \S 1.  

\bigskip

The structure of this text is as follows. We state in \S \ref{sec1} our main result and explain the connection with a conjecture of the second author. We gather in \S \ref{sec2} the main tools in the proof, namely the multiplicity and interpolation estimates we rely on, and also a counting lemma which provides an asymptotic estimate for the cardinality of the image of $\GS $ in sub-quotients of $G$. Then we construct in \S \ref{sec3} the chain of algebraic subgroups $\chaine$ and study its properties. This section might be of independent interest, and is logically independent from the previous ones. Finally in \S \ref{sec4} we prove our main result and gather in \S \ref{secgen} some remarks and comments on possible generalizations.

\bigskip

\noindent
{\it Acknowledgments}  It is a pleasure to thank the Universit\'e de Paris-Sud, Orsay
for receiving the second author during January and February, 2012, providing an 
opportunity to start this work.  The second author would also like to thank Imperial
College which provided a pleasant environment in which to continue 
 work on this article. The first author is partially supported by Agence Nationale de la Recherche (project HAMOT, ref. ANR 2010 BLAN-0115-01), and both would like to warmly
 thank Michel Waldschmidt for  his long standing
 encouragement and his stimulating questions.

\section{Statement of the results} \label{sec1}

Throughout this section we  let $G$ be a connected commutative algebraic group embedded in  projective space $\P^N$.  Suppose  $\gamul\in G$ and let   $\Gamma$ denote  the subgroup generated by $\gamul$.
Let $\Sulgeq \geq 1$ be real numbers, and recall that 
$\GS$ is the set of all points $n_1\gamma_1+\ldots+n_l\gamma_l$ with integers $n_j$ such that $|n_j| < S_j$; here  $\Ssoul$ denotes the tuple $(\Sul)$.

Using  this data we shall construct in  \S \ref{sec3} a   chain of algebraic subgroups $\{0\} = H_0 \subsetneq H_1  \subsetneq \ldots   \subsetneq H_r = G$, with $r\geq 1$.

\bigskip

We let $\Gamma_j$ denote  the subgroup generated by $\gamma_1$, \ldots, $\gamma_j$, setting
 $\Gamma_0 = \{0\}$, and we put
$$\sss_i = \left(  \prod_{j=1}^l  S_j^{\rk \Big(\frac{\Gamma_j\cap H_{i+1}}{\Gamma_j\cap H_i}\Big) - \rk\Big(\frac{\Gamma_{j-1}\cap H_{i+1}}{\Gamma_{j-1}\cap H_i}\Big)} \right)^{1/(\dimhipu-\dimhi)} $$
for any $i \in \zerormu$. Then we shall prove that
\begin{equation} \label{eqordresi}
1 \leq \sss_{r-1} <  \sss_{r-2} < \ldots < \sss_1 < \sss_0,
\end{equation}
as an immediate consequence of Proposition \ref{propchaineiii} and Eq. \eqref{eqgensi} in \S \ref{sec3}.

\subsection{The main result}

Our main result is twofold. The first one is proved using  interpolation  estimates, whereas the second one is based on multiplicity estimates.

\begin{Th} \label{th11}
There exists a positive constant $\cthuu$, depending only on $G \hookrightarrow \P^N$,   $\gamul$ but not on $\Sulgeq$, such that
$$\base \subset \GS + H_i$$
 for any positive integers $D$, $T$ such that  $D> \cthuu \sss_i T$ with $i \in \zerormu$.
 \end{Th}

\begin{Th} \label{th12}
For any $\eps  $   with $0 < \eps < 1$  there exists a positive constant $\cthud$, depending only on $G \hookrightarrow \P^N$, $\eps$, $\gamul$ but not on $\Sulgeq$, such that
$$\Gamma((1-\eps)\Ssoul) +H_i \subset \base  $$
 for any positive integers $D$, $T$ such that  $D< \cthud^{-1} \sss_{i-1} T$ with $i \in \unr$.
 \end{Th}

Theorem \ref{th12} is closely related to Lemma 1.5.3 in \cite{NakamayeLang}.  This latter result
assumes that $S_1 = S_2 = \ldots = S_l$ and it
only treats the case $i = 1$. It is stated for $H_1$ alone rather than $\Gamma((1-\eps)\Ssoul) + H_1$
but it applies to these translates of $H_1$.  Subgroups closely related to the
sequence $H_2,\ldots,H_r$  appear in Conjecture 1.1.9 of \cite{NakamayeLang}.  The techniques
of \cite{NakamayeLang} are completely different from the present paper.  In particular all 
constructions take place on $X$: no embeddings or quotient maps are used.  The end result
is that the results of \cite{NakamayeLang} are quantitatively stronger (the constants are
sharp in the same way as those of Philippon's multiplicity estimates)  but they apply in
very few cases. 
 
 \bigskip
 
 The cases where $D$ is very small or very large in comparison with $T$ and the $\sss_i$ will be dealt with in \S \ref{subsec21}. If $D <  \cthud^{-1} \sss_{r-1} T$ then Theorem \ref{th12} asserts that $\base = G$; this is a multiplicity estimate, stated below as Proposition \ref{proppphz}. In a ``dual'' way, if  $D >  \cthuu \sss_{0} T$ then Theorem \ref{th11} means that $\base =  \GS $ since the inclusion $\GS \subset \base$ holds trivially. We shall derive this result, stated as Proposition \ref{propinterp}, from an interpolation estimate, namely Proposition \ref{propinterpFN}. 

When $r = 1$, we do not prove anything more -- we could probably refine our result in this case, to make the constants $\cthuu$ and $\cthud$ explicit, but we are not able to do it in general (see \S \ref{secgen}). When $r \geq 2$, our proof procedes by applying these results in sub-quotients of $G$ coming from the algebraic subgroups $H_i$.

If $r\geq 2$ and 
$$\cthuu \sss_i T < D < \cthud^{-1}  \sss_{i-1}  T $$
 for some $i\in\unrmu$,  which happens for some integers $D$  provided $\Card \,  \GS$  is sufficiently large in terms of $G \hookrightarrow \P^N$, $\eps$, $\gamul$, then
$$\Gamma((1-\eps)\Ssoul)+H_i\subset \base \subset \GS + H_i.$$
Therefore Theorem \ref{thintro} follows from Theorems \ref{th11} and \ref{th12}, since when $S_1=\ldots=S_l=S$ we have
$$\sss_i =  S^{\mu_i} \mbox{ with } \mu_i =    \frac{\rk(\Gamma\cap H_{i+1}) - \rk(\Gamma\cap H_{i})}{\dim H_{i+1} - \dim H_i} .$$

\begin{Rem}  We shall prove in Lemma \ref{lemab} below (\S \ref{subsec22}) that $\sss_i^{\dimhipu-\dimhi}$ is equal to the cardinality of $(\GS \cap H_{i+1} )  \bmod H_i$, up to a multiplicative constant depending only on $\gamul$.  Therefore $\sss_i$ might be replaced by the $(\dimhipu-\dimhi)$-th root of this cardinality in Theorems \ref{th11} and \ref{th12}  up to changing the values of the constants $\cthuu$ and $\cthud$.  The assumption  $D< \cthud \sss_{i-1} T$ in Theorem  \ref{th12} is the one needed to apply a multiplicity estimate in $H_{i }/H_{i-1}$ in order to guarantee that no non-zero polynomial of degree $D$ on $H_{i}/H_{i-1}$ vanishes to order at least $T$ at each point of  $(\GS \cap H_{i } )  \bmod H_{i-1}$. Of course $\cthud$ should take here a suitable value in terms of a projective embedding of $H_{i}/H_{i-1}$. The same remarks apply to the assumption  $D> \cthuu \sss_{i } T$ in Theorem  \ref{th11} needed to apply an interpolation estimate (or Proposition \ref{propinterp} below) on
  $(\GS \cap H_{i+1} )  \bmod H_i$ in the algebraic group $H_{i+1}/H_i$ (see \S \ref{subsec21}).
\end{Rem}

\subsection{Connection to a conjecture of the second author} \label{subsecconj}

Following \cite{NakamayeLang} we let 
$$\alpha_j = \sup\{\alpha\in\Q,\,\dim \basek < j \mbox{ for any $k$ sufficiently large}\}$$
where $j\in\unn$ and $n = \dim G$.  Theorems \ref{th11} and \ref{th12}, applied with $\eps=1/2$,   yield 
$$\cthuu^{-1} \sss_i^{-1} \leq\alpha_j\leq \cthud  \sss_i^{-1}$$
where $i\in\zerormu$ is  chosen so that   $\dimhi < j \leq \dimhipu$.  Consequently 
$$ \cthuu^{-n}  \prod_{i=0}^{r-1} \sss_i^{-(\dimhipu-\dimhi)}  \leq
\prod_{j=1}^n \alpha_j \leq  \cthud^n  \prod_{i=0}^{r-1} \sss_i^{-(\dimhipu-\dimhi)}.$$
Thus
$$\cthuu^{-n} \Big[ \prod_{j=1}^l S_j^{\rk\, \Gamma_j - \rk \, \Gamma_{j-1}}\Big]^{-1}  \leq
\prod_{j=1}^n \alpha_j \leq   \cthud^n \Big[ \prod_{j=1}^l S_j^{\rk \,\Gamma_j - \rk\, \Gamma_{j-1}}\Big]^{-1}.$$
Using Lemma \ref{lemab} below with $H' = G$ and $H''=\{0\}$ we obtain a positive constant $\cconj$, depending only on $G \hookrightarrow \P^N$ and $\gamul$, such that
$$\cconj^{-1} \leq \Big(\Card \,  \GS\Big) \prod_{j=1}^n \alpha_j \leq \cconj.$$
Of course the important point here is that $\cconj$ does not depend on $\Sul$. In parallel to Conjecture 1.1.4 of \cite{NakamayeLang},   it seems natural to ask whether
$$\frac{\deg_{\oo(1)}(\Gbar)}{n!} \leq \Big(\Card \, \GS\Big) \prod_{j=1}^n \alpha_j \leq \deg_{\oo(1)}(\Gbar),$$
where $n = \dim G$ and $\Gbar$ is the Zariski closure of $G \hookrightarrow \P^N$.  The upper bound can be proved using intersection theory and the definition
 of the $\alpha_i$,  as in  \S 1.2 of  \cite{NakamayeLang}.

\section{Prerequisites} \label{sec2}

In this section we state the interpolation and multiplicity estimates we rely on and apply them to the extremal cases $i=0$ (in Theorem \ref{th11}) and $i=r$ (in Theorem \ref{th12}). Then we state and prove in \S \ref{subsec22} a lemma that provides an asymptotic estimate for the cardinality of the image of $\GS$ in sub-quotients of $G$.

\subsection{Interpolation and Multiplicity Estimates} \label{subsec21}

We shall use the following notation: given a finite subset $\Omega$ of a commutative algebraic group $G$ and a positive integer $n$, we let 
$\Omega[n]$ denote the set of all sums $\omega_1 + \ldots + \omega_n$ where $\omega_1$, \ldots, $\omega_n$ are (not necessarily distinct) elements of $\Omega$. We denote by  $\Omega\{ n\}  $ the set $ \Omega [n]  -  \Omega [n] $, that is the set of all elements $x-y$ with $x,y\in   \Omega [n] $.

\bigskip

The following is a  weak form of the  multiplicity estimate, Theorem~2.1, from \cite{PPhzeros}.

\begin{Prop} \label{proppphz}
Let $G$ be a connected commutative algebraic group, embedded in  projective space $\P^M$. Then there is a positive constant $\cpphz$, depending only on $G$ and on this embedding, with the following property.
Let $\Omega$ be a finite subset of $G$, and suppose $D$, $T$ are
 positive integers such that, for  every \cas $H \subsetneq G$, 
\begin{equation} \label{eqhypproppphz}
\Card(\Omega \bmod H) \, T^{\dim(G/H)} > \cpphz D ^{\dim(G/H)} .
\end{equation}
Then no non-zero $P \in\cdegD$ vanishes to order at least $T$ at every point of $\Omega[\dim G]$. In other words,
$$\baseom = G.$$
\end{Prop}

\bigskip

We shall deduce the statement ``dual'' to Proposition \ref{proppphz}, namely Proposition \ref{propinterp}, from the following interpolation estimate (which is Corollary 1.2 of \cite{FNinterp}).

\begin{Prop} \label{propinterpFN}
Let $G$ be a connected commutative algebraic group, embedded in   projective space $\P^M$. Then there is a positive constant $\clemfn$, 
depending only on $G$ and on this embedding, with the following property.
Let $\Omega$ be a finite subset of $G$, and suppose $D$, $T$ are positive integers such that, for any translate $x+H$ of a non-zero \cas $H  $ of $G$, 
$$  \Card \Big((\Omega\cap (x+H))[\dim(H)] \Big)  \,  T^{ \dim(H)}  <  \clemfn  D^{ \dim(H)}  .$$
 Then the evaluation map
$$
\cdegD = H^0(\Gbar ,\oo(D)) \to  H^0\left(\Gbar , \oo(D) \otimes \oplus_{\omega \in \Omega}\oo_\Gbar/m_\omega^{T}
 \right)
$$
is surjective, where $\Gbar$ is the Zariski closure of $G$ in $\P^M$ and $m_\omega \subset \oo_\Gbar$  is the maximal ideal sheaf corresponding to the point $\omega$.
\end{Prop}

This result is essentially as precise as Philippon's multiplicity  estimate (namely   Theorem 2.1 of \cite{PPhzeros}), and even slightly more. A less precise estimate (in the style of \cite{Masser} or \cite{SFinterp}) would not be sufficient to deduce the following result, which we shall use later in this text.

\begin{Prop} \label{propinterp}
Let $G$ be a connected commutative algebraic group, embedded in   projective space $\P^M$. Then there is a positive constant $ \cinterp $, depending only on $G$ and on this embedding, with the following property.
Let $\Omega$ be a finite subset of $G$, and suppose $D$, $T$ are positive integers such that, for every non-zero \cas $H  $ of $G$, 
\begin{equation} \label{eqhyppropinterp}
\Card(\Omega\{n\} \cap H)  \, T^{\dim H } < \cinterp D ^{\dim H }  
\end{equation}
where $n = \dim G$. Then 
$$\baseom = \Omega.$$
\end{Prop}

It should be noticed that only algebraic subgroups $H$ appear in this result, whereas translates are needed in Proposition \ref{propinterpFN}. This is due to the fact that $\Omega\{n\} $ (i.e.,  the set   of all elements $x-y$ with $x,y\in   \Omega [n] $) is used instead of $\Omega[n]$.

\bigskip

\Dem of Proposition \ref{propinterp} :  The inclusion $\Omega \subset \baseom$ holds trivially. Let $g\in G \setminus \Omega$ and put $\Omega' = \Omega\cup\{g\}$. Let $H' = x+H$ be any translate  of a non-zero \cas $H  $ of $G$. Then $(\Omega'\cap H')[\dim H] \subset \cup_{i=0}^{\dim H} \calE_i$ where
$$\cale_i = \{ i g + \gamma, \, \gamma \in \Omega[\dim H - i]\} \cap H'.$$
If $\cale_i \neq \emptyset$, substracting a fixed element of $\cale_i$ yields  an injective map 
$$\cale_i \to \Omega\{\dim H -i\} \cap H  \subset  \Omega\{n\} \cap H ,$$
so that $\Card \,  \cale_i \leq \Card (  \Omega\{n\} \cap H  )$, and this inequality holds also if $\cale_i = \emptyset$. Therefore we have
$$ \Card \Big((\Omega ' \cap  H' )[\dim(H)] \Big)     \,  T^{ \dim(H)}   \leq (n+1)  \Card \Big(\Omega\{n\} \cap H  \Big) \,  T^{ \dim(H)} <  (n+1)  \cinterp D ^{\dim H }  .$$
Choosing $\cinterp =  \clemfn / (n+1)$, Proposition \ref{propinterpFN} provides $P \in\cdegD$ which vanishes to order at least $T$ at each point of $\Omega$ and does not vanish at $g$. This proves that $g \not\in\baseom$, and concludes the proof of Proposition \ref{propinterp}.

\subsection{A Counting Lemma} \label{subsec22}
 
The following lemma is very useful for estimating the number of points of $\GS$ in sub-quotients of $G$. The fundamental idea is that $\Sul$ will be assumed to be sufficiently large, in terms of $\gamul$, so that this number of points can be estimated asymptotically in terms of ranks of $\Z$-modules. Recall that $\Gamma_j$ denotes the  $\Z$-module  generated by $\gamma_1$, \ldots, $\gamma_j$, with $\Gamma_0 = \{0\}$.

\begin{Lemme} \label{lemab}
Let $H'$, $H''$ be algebraic subgroups of a commutative algebraic group $G$, such that $H'' \subset H'$. Let $\gamul\in G$ and 
let $\Gamma$ be the subgroup generated by $\gamul$. Then there exist positive constants $\caba$ and $\cabb$ with the following properties:
\begin{itemize}
\item $\caba$ depends only on $\gamul$ and on $H'$ (but not on $H''$).
\item $\cabb$ depends only on $\gamul$ and on $H''$ (but not on $H'$).
\item For any real numbers $\Sulgeq \geq 1$ we have
$$\caba \Sgothhprhsec < \Card (\GS \cap H' \bmod  H'') < \cabb \Sgothhprhsec$$
where
$$\Sgothhprhsec = \prod_{j=1}^l S_j^{\rk \Big(\frac{\Gamma_j\cap H'}{\Gamma_j\cap H''}\Big) - \rk\Big(\frac{\Gamma_{j-1}\cap H'}{\Gamma_{j-1}\cap H''}\Big)}.$$
\end{itemize}
\end{Lemme}

In the special case $H'' = \{0\}$, Lemma \ref{lemab} reduces to Lemma 1.5 of \cite{SFinterp}, except that in \cite{SFinterp} the constant $\cabb$ may depend on $H'$.

We did not try to make explicit the constants $\caba$ and  $\cabb$ since it is not needed 
in our application. 
However it is critical that $\caba$ does not depend on $H''$ and that $\cabb$ does not
depend on $H'$.

To illustrate this situation, let us consider the case where $l=1$ and $\gamma_1$ is not torsion in $G$. Let $H$ be a \cas of $G$ which contains $N\gamma_1$ for some $N\geq 1$, but not $k\gamma_1$ with $1\leq k \leq N-1$.  Then $\Card ( \Gamma(S_1) \cap H) = 2M+1$, where $M \geq 0$ is the largest integer such that $MN < S_1$, and $\Card (  \Gamma(S_1) \bmod H) = N$  if $S_1> N $. Taking $H' = H$ and $H'' = \{0\}$ we see that $\caba$ has to depend on $H'$, since $N$ may take arbitrarily large values in terms of $\gamul$. In the same way, taking $H'=G$ and $H''=H$ shows that $\cabb$ has to depend on $H''$.
 
Using Lemma \ref{lemchainev} and the notation of \S \ref{sec1},   Lemma \ref{lemab}  proves  that $\sss_i^{\dimhipu-\dimhi}$ is equal to the cardinality of $(\GS \cap H_{i+1} )  \bmod H_i$, up to a multiplicative constant depending only on $\gamul$.

\begin{Rem} \label{remab}
An immediate consequence of Lemma \ref{lemab} is that for any $x\in G$ 
$$  \Card (\GS \cap (x+H') \bmod  H'') < \cabc \Sgothhprhsec$$
where $\cabc$ depends only on $\gamul$ and on $H''$ but not on $H'$ or on $x$. Indeed, subtracting a fixed element of $\GS \cap (x+H')$ yields an injective map
$\GS \cap (x+H') \to \Gamma(2\Ssoul) \cap  H' $. 
\end{Rem}

\begin{Rem} \label{remabd} For  any $\lambda\geq 1$, applying Lemma \ref{lemab} with $\lambda \Ssoul = (\lambda S_1,\ldots,\lambda S_l)$ yields
\begin{equation} \label{eqabd}
\caba \lambda^{\rk  \Big(\frac{\Gamma\cap H'}{\Gamma \cap H''}\Big)} \Sgothhprhsec < \Card (\Gamma(\lambda \Ssoul) \cap H' \bmod  H'') < \cabb  \lambda^{\rk  \Big(\frac{\Gamma\cap H'}{\Gamma \cap H''}\Big)} \Sgothhprhsec.
\end{equation}
This will be used several times in the proof of Lemma \ref{lemab}, without explicit reference. Moreover the first inequality in Eq. \eqref{eqabd} holds for any $\lambda > 0$.
\end{Rem}

\bigskip

\Dem of Lemma \ref{lemab}: Since the result is trivial when $l=0$, we may assume by induction that it holds for $\glmu$. Notice that the value of $\Sgothhprhsec$ relative to $\GS$ is the same as the one relative to $\glmu(S_1,\ldots,S_{l-1})$ if $\rk \Big(\frac{\Gamma \cap H'}{\Gamma \cap H''}\Big) =  \rk\Big(\frac{\glmu \cap H'}{\glmu \cap H''}\Big)$, and it is $S_l$ times bigger otherwise.

The lower bound on $ \Card (\GS \cap H' \bmod  H'')$ follows at once from the inclusion \linebreak $\glmu(S_1,\ldots,S_{l-1}) \subset \GS$ if $\rk \Big(\frac{\Gamma \cap H'}{\Gamma \cap H''}\Big) =  \rk\Big(\frac{\glmu \cap H'}{\glmu \cap H''}\Big)$. Otherwise we have $\rk( \Gamma \cap H' ) =  1+ \rk( \glmu \cap H' ) $ and $\rk( \Gamma \cap H'' ) =   \rk( \glmu \cap H'' ) $. In this case, there exist $m_1,\ldots,m_l\in\Z$, with $m_l\geq 1$, such that $\gti  = m_1\gamma_1+\ldots+m_l\gamma_l$ belongs to $\Gamma\cap H'$ and has infinite order in $\frac{ \Gamma\cap H' }{\glmu \cap H' }$. Letting $M = \max(|m_1|, \ldots, |m_l|)$, the elements $\gamma_0 + n\gti$, where $|n| < S_l / 2M$ and $\gamma_0$ ranges through a system of representatives of $\glmu(S_1/2,\ldots,S_{l-1}/2)  \cap H' \bmod H''$, belong to $\GS\cap H'$ since $S_l \leq S_i$ for any $ i \in \unlmu$. If two of them are equal modulo $H''$, say $\gamma_0 + n\gti \in \gamma'_0 + n'\gti + H''$ with $(\gamma_0 , n) \neq (\gamma'_0 , n') $, then $(n-n')\gti + (\gamma_0-\gamma'_0)\in\Gamma \cap H''$. Now $\glmu \cap H''$ has finite index, say $N$, in $\Gamma\cap H''$ so that $N(n-n')\gti + N(\gamma_0-\gamma'_0)\in\glmu \cap H'' \subset \glmu \cap H'$ and $N(n-n')\gti \in\glmu \cap H'$. Since the image of $\gti$ in $\frac{ \Gamma\cap H' }{\glmu \cap H' }$ has infinite order, we have $n=n'$ and $\gamma_0  \in \gamma'_0  + H''$ which is a contradiction. Therefore the elements given above are pairwise distinct, concluding the proof of the  lower bound on $ \Card (\GS \cap H' \bmod  H'')$.

To prove the upper bound, we distinguish between three cases.

$(a)$ If  $\rk( \Gamma \cap H' ) =   \rk( \glmu \cap H' ) $, the upper bound holds trivially if we 
 also have $\GS\cap H' = \glmu(S_1,\ldots,S_{l-1})\cap H'$. Otherwise 
 there exist $m_1,\ldots,m_l\in\Z$, with $m_l\geq 1$, such that $\gti  = m_1\gamma_1+\ldots+m_l\gamma_l$ belongs to $\Gamma\cap H'$.  Letting $N$ denote the index of $ \glmu \cap H' $ in $ \Gamma \cap H' $, we have $N\gti = Nm_1\gamma_1+\ldots+Nm_l\gamma_l \in  \glmu \cap H' $  so that $Nm_l\gamma_l \in  \glmu $. Therefore the image of $\gamma_l$ in $\Gamma/\glmu$ has finite order: let $\omega$ denote this order, which depends only on $\gamul$. There exist $r_1,\ldots,r_l\in\Z$ such that $r_1\gamma_1+\ldots+r_l\gamma_l = 0$ where $r_l = \omega\geq 1$. Letting 
  $R = \max(|r_1|, \ldots, |r_l|)$, we have
  $$\GS \subset \cup_{n=0}^{r_l-1}n\gamma_l + \glmu((R+1)S_1,\ldots,(R+1)S_{l-1})$$
since $S_l\leq S_i$ for any $i\in\unlmu$. The upper bound follows at once.

$(b)$ If $\rk( \Gamma \cap H' ) =  1+ \rk( \glmu \cap H' ) $ and  $\rk( \Gamma \cap H'' ) =   \rk( \glmu \cap H'' ) $,  there exist $m_1,\ldots,m_l\in\Z$ such that $\gti  = m_1\gamma_1+\ldots+m_l\gamma_l\in \Gamma\cap H'$ and $m_l \geq 1$; we choose these integers with the least possible value of $m_l$. Then for any $\gamma = n_1\gamma_1+\ldots+n_l\gamma_l\in\GS\cap H'$, $n_l$ is a multiple of $m_l$ and we have $\gamma - r\gti \in\glmu\cap H'$ where $r = n_l / m_l$ is such that $|r|< S_l$. Letting $\gamma' = \gti - m_l\gamma_l  =  m_1\gamma_1+\ldots+m_{l-1}\gamma_{l-1}$ we obtain
$$\GS \subset \cup_{r=-S_l}^{S_l} rm_l\gamma_l + \Big[  \glmu(S_1,\ldots,S_{l-1}) \cap (r\gamma'+H')\Big].$$
Using Remark \ref{remab} this concludes the proof of the upper bound in this case.

$(c)$ If $\rk( \Gamma \cap H' ) =  1+ \rk( \glmu \cap H' ) $ and  $\rk( \Gamma \cap H'' ) =  1+ \rk( \glmu \cap H'' ) $,   there exist $m_1,\ldots,m_l\in\Z$ such that $\gti  = m_1\gamma_1+\ldots+m_l\gamma_l\in \Gamma\cap H''$ and $m_l \geq 1$. Let $\gamma = n_1\gamma_1+\ldots+n_l\gamma_l\in\GS\cap H'$, and let $q,  r \in\Z$ be such that $n_l = qm_l+r$ with $|r_l|<m_l$ and $|q|\leq \frac{|n_l|}{m_l}<S_l$. Then we have $\gamma-q\gti = (n_1-qm_1)\gamma_1+\ldots+  (n_{l-1}-qm_{l-1})\gamma_{l-1}+r\gamma_l$ so that, letting $M =  \max(|m_1|, \ldots, |m_l|)$,
$$\GS \cap H' \bmod  H'' \subset \cup_{r=-M}^M r\gamma_l + \Big[ \glmu((M+1)S_1,\ldots,(M+1)S_{l-1})\cap(-r\gamma_l+H')\Big] \bmod H''.$$
Using Remark \ref{remab} this concludes the proof of  Lemma \ref{lemab}.

\section{A Chain of Algebraic Subgroups} \label{sec3}

Throughout this section we fix a connected  commutative algebraic group $G$, real  numbers $\Sul$ and elements $\gamul \in G$; we assume that $\Sulgeq \geq 1$. With this data we associate in \S \ref{subsec31} a chain of \cass $\chaine$ of $G$. We study its properties throughout this section, with a special emphasis on its connection to the distribution of $\GS$ with respect to algebraic subgroups of $G$ (\S \ref{subsec32}), and on the case where $S_1 = \ldots = S_l$ as in the introduction (\S \ref{subsec33}).

\subsection{Construction and First Properties} \label{subsec31}

For any \cas $K$ of $G$ we let 
$$\fis(K) = \sum_{j=1}^l \Big(  \rk(\Gamma_j\cap K)- \rk(\Gamma_{j-1}\cap K)\Big) \log S_j,$$
where $\Gamma_j$ is the subgroup of $\Gamma$ generated by $\gamma_1,\ldots,\gamma_j$,
 $\Gamma_0 = \{0\}$, and $\Ssoul = (\Sul)$.
With this definition, $\Card(\GS \cap K) $ is essentially equal to $\exp \fis(K)$ by
 Lemma \ref{lemab} above, with $H'' = \{0\}$, so that $\SgothKzero = \exp \fis(K)$. We refer to  \S 3 of  \cite{SFinterp} for a related construction.

\bigskip

In the special case where $S_1=\ldots=S_l=S$, we have $\fis(K) = \rk(\Gamma \cap K)\log S $. The starting point of our construction is the existence  \cite{NakamayeLang}, in this case, of a maximal element $H_1$ with respect to inclusion among the non-zero \cass $H$ such that  $\muet(\Gamma,G) =\frac{ \rk(\Gamma\cap H)}{  \dim H}$. Reapplying this  construction 
in $G/H_1$ with $\Gamma\bmod H_1 = (\Gamma+H_1)/H_1$ yields a maximal \cas $H_2/H_1$ of $G/H_1$, with $H_1 \subsetneq H_2$. Repeating this argument leads to a chain of algebraic subgroups of $G$, which we construct now in the general case where $\Sul$ are not assumed to be equal.

\begin{Prop} \label{propchainez}
There exists a unique chain $\{0\} = H_0 \subsetneq H_1  \subsetneq \ldots   \subsetneq H_r = G$ of \cass  of $G$, with $r\geq 1$, such that:
\begin{itemize}
\item For any $i\in\zerormu$ and any \cas $K$
such that $\dim K > \dim H_i$, we have
\begin{equation} \label{eqchaine1}
\frac{\fis(K) - \fis(H_i)}{\dim K - \dim H_i} \leq \frac{\fis(H_{i+1}) - \fis(H_i)}{\dim H_{i+1} - \dim H_i}.
\end{equation}
\item If equality holds in  Eq. \eqref{eqchaine1} then $H_i \subset K \subset H_{i+1}$.
\end{itemize}
\end{Prop}

\begin{Rem} \label{remchainei} In the proof of Proposition \ref{propchainez} we shall prove actually a stronger property of these subgroups, namely that 
for any $i\in\zerormu$ and any  \cas   $K$ we have 
\begin{eqnarray} 
 &&[\dimhipu - \dimhi] \fis(K) - [\fishipu - \fishi] \dim K \nonumber \\
&+& \fishipu \dimhi - \fishi \dimhipu \leq 0 \label{eqchaine2}.
\end{eqnarray}
If equality holds then $H_i \subset K \subset H_{i+1}$.
This  inequality can be also be written as
$$ [ \fis(K) - \fis(H_i) ] [ \dim H_{i+1} - \dim H_i ] \leq [ \dim K - \dim H_i ] [ \fis(H_{i+1}) - \fis(H_i) ] .$$
If  $\dim K > \dim H_i$  it is equivalent to Eq. \eqref{eqchaine1}. If $\dim K < \dim H_i$ it yields
\begin{equation} \label{eqchainenv}
\frac{\fis(H_i) - \fis(K)}{\dim H_i - \dim K}  > \frac{\fis(H_{i+1}) - \fis(H_i)}{\dim H_{i+1} - \dim H_i}.
\end{equation}
\end{Rem}

In the case where $S_1=\ldots=S_l$, reasoning as in   \cite{NakamayeLang} one can prove the existence of a minimal element $H_{r-1}$ with respect to inclusion among the  \cass $H \subsetneq G$ such that  $\mu(\Gamma,G) =\frac{ \rk(\Gamma) - \rk(\Gamma\cap H)}{   \dim G -  \dim H}$. Applying this property again in $H_{r-1}$ with $\Gamma\cap H_{r-1}$ provides $H_{r-2} \subsetneq H_{r-1}$. The following immediate consequence of Eq. \eqref{eqchaine2} in Remark  \ref{remchainei} asserts that the chain of \cass of $G$ constructed by iterating this process (and generalizing it to allow $\Sul$ not to be equal)  is the same as above.

\begin{Prop} \label{propchaineii}
For any $i\in\zerormu$ and any  \cas   $K$  such that $\dim K  < \dimhipu$ we have
\begin{equation} \label{eqchaine3}
  \frac{\fis(H_{i+1}) - \fis(K)}{\dim H_{i+1} - \dim K } \geq  \frac{\fis(H_{i+1}) - \fis(H_i)}{\dim H_{i+1} - \dim H_i}.
\end{equation}
Moreover if equality holds then $H_i \subset K \subset H_{i+1}$.
\end{Prop}

Assuming again $S_1=\ldots=S_l$, $H_1$ is maximal such that $\frac{ \rk(\Gamma\cap H_1)}{  \dim H_1} = \muet(\Gamma,G) =\max_{H\neq 0} \frac{ \rk(\Gamma\cap H)}{  \dim H}$. In particular we have $\frac{ \rk(\Gamma\cap H_1)}{  \dim H_1}  > \frac{ \rk(\Gamma\cap H_2)}{  \dim H_2} $  since $H_1\subsetneq H_2$. Now  $\frac{ \rk(\Gamma\cap H_2)}{  \dim H_2}$ lies between  $\frac{ \rk(\Gamma\cap H_1)}{  \dim H_1}$ and  $\frac{ \rk(\Gamma\cap H_2) -  \rk(\Gamma\cap H_1)}{  \dim H_2 -  \dim H_1}$ because its numerator is the sum of both numerators and the same property holds for the denominators, so that  $\frac{ \rk(\Gamma\cap H_2) -  \rk(\Gamma\cap H_1)}{  \dim H_2 -  \dim H_1} < \frac{ \rk(\Gamma\cap H_1)}{  \dim H_1}$. Generalizing this result to all subgroups $H_i$ and removing the assumption $S_1=\ldots=S_l$, we obtain the following.

\begin{Prop} \label{propchaineiii}
For any $i\in\unrmu$ we have
$$ \frac{\fishipu - \fishi}{\dimhipu - \dimhi} < \frac{\fishi  - \fishimu}{\dimhi  - \dimhimu}.$$
\end{Prop}

This proposition  follows immediately from Eq. \eqref{eqchainenv} by taking $K=H_{i-1}$. It    is the key point in the proof of Eqns. \eqref{eqmudecr} and \eqref{eqordresi} above.

\begin{Rem} \label{rempolygone}
With each \cas $K$ of $G$ we may associate the point $M_K = (\dim K , \fis(K))\in\R^2$. Then  our construction  yields a convex polygon $M_{H_0}M_{H_1}\ldots M_{H_r}N$, where $N = (\dim G, 0)$. For any $K$ the point $M_K$ is either inside this polygon or on an edge; if it lies on the segment $[M_{H_i} M_{H_{i+1}}]$ with $0 \leq i \leq r-1$ then $H_i \subset K \subset H_{i+1}$. Indeed Eq. \eqref{eqchaine1} means that the line $(M_{H_i} M_{K})$ has slope less than or equal to that of $(M_{H_i} M_{H_{i+1}})$, if $\dim K > \dimhi$. This means that $M_K$ is below the line  $(M_{H_i} M_{H_{i+1}})$, which is expressed by Eq. \eqref{eqchaine2}. An equivalent statement, if $\dim K < \dimhipu$, is provided by Proposition \ref{propchaineii}. Namely, the slope of  $(M_K M_{H_{i+1}})$ is less than or equal to that of $(M_{H_i} M_{H_{i+1}})$. Lastly the slope of $(M_{H_i} M_{H_{i+1}})$ is a decreasing function of $i$, as Proposition \ref{propchaineiii} states.
\end{Rem}

\Dem of Proposition \ref{propchainez}   and Remark  \ref{remchainei}:  To begin with, we notice hat 
$$\fis(K) = \sum_{j=1}^l \rk(\Gamma_j\cap K) \log(S_j/S_{j+1})$$
for any \cas $K$ of $G$: here we let $S_{l+1}=1$. For any $j$ and any  \cass  $K$, $K'$ of $G$ 
 we have  $\rk(\Gamma_j\cap(K\cap K')) + \rk(\Gamma_j\cap(K+ K')) \geq \rk(\Gamma_j\cap K) +  \rk(\Gamma_j\cap K')$ so that
\begin{equation} \label{eqsommeinter1}
\fis(K \cap K')+\fis(K+K')\geq \fis(K) + \fis(K')
\end{equation}
since $\log(S_j/S_{j+1})\geq 0$ for any $j$. We shall also use the fact that 
\begin{equation} \label{eqsommeinter2}
\dim(K \cap K')+\dim(K+K') = \dim(K) + \dim(K').
\end{equation}
Now let us construct $H_i$ and prove the results at the same time, by induction on $i$. If the algebraic subgroups $H_0$, \ldots, $H_i$ satisfy the desired properties with $i\geq 0$ and $H_i \neq G$, we define $H_{i+1}$ to be a \cas of $G$ of dimension greater than $\dim H_i$ for which $\frac{\fishipu-\fishi}{\dimhipu-\dimhi}$ is maximal. If there are several \cass $K$ of $G$ with $\dim K > \dimhi$ for which $\frac{\fis(K) - \fishi}{\dim(K) - \dimhi}$ is equal to this maximal value, then we choose $H_{i+1}$ with maximal dimension among them. In this way Eq. \eqref{eqchaine1} holds for any $K$ such that  $\dim K > \dimhi$.  Moreover the \cas  $H_{i+1}$ constructed in this way is unique: if a \cas $H'_{i+1}$ satifies $\frac{\fis ( H'_{i+1} ) -\fishi}{\dim H'_{i+1}-\dimhi} = \frac{\fishipu-\fishi}{\dimhipu-\dimhi}$ and $\dim H'_{i+1} = \dim H_{i+1}$ then $ H'_{i+1} = H_{i+1}$; this follows from the case of equality in Eq.   \eqref{eqchaine1}, which will be proved below. As a consequence,   the chain $\chaine$ is unique.   It should also be emphasized that $H_{i+1}$ is constructed in terms of the function $\fis$ only; accordingly it depends only on $\Gamma_1$, \ldots, $\Gamma_l$, not really on $\gamul$. 

\bigskip

Now let $\chi(K)$ denote the left handside of Eq. \eqref{eqchaine2}; notice that 
\begin{equation} \label{eq8bis}
\chi(H_i) = \chi(H_{i+1})=0.
\end{equation}
 Actually if we associate 
with each \cas $K$ of $G$   the point $M_K = (\dim K , \fis(K))$ as in Remark \ref{rempolygone} above, then $\chi(K) = 0$ means that $M_K$ lies on the line   $(M_{H_i} M_{H_{i+1}})$.

By definition of $H_{i+1}$ we have
\begin{equation} \label{eqchi1}
\left\{\begin{array}{l} \chi(K) \leq 0 \mbox{ if } \dim K > \dimhi\\
\mbox{if equality holds then } \dim K \leq \dimhipu. \end{array}\right.
\end{equation}
To conclude the proof of  Eq. \eqref{eqchaine2} for any $K$, let us prove also that
\begin{equation} \label{eqchi2}
\left\{\begin{array}{l} \chi(K) \leq 0 \mbox{ if } \dim K \leq \dimhi\\
\mbox{if equality holds then }K=H_i. \end{array}\right.
\end{equation}
If $i=0$ then this is a triviality. If $i\geq 1$, $\dim K = \dimhi$ and $K \neq H_i$ then $\chi(K) = (\dimhipu-\dimhi)(\fis(K)-\fis(H_i))<0$ using Eq. \eqref{eqchaine1} with $i-1$. If $i\geq 1$ 
and  $\dim K < \dimhi$, notice that Eq. \eqref{eqchaine2} with $i-1$ reads
$$\frac{\fishi -\fis(K)}{\dimhi - \dim K} \geq \frac{\fishi -\fishimu}{\dimhi - \dimhimu}.$$
Combining this inequality with
Proposition \ref{propchaineiii} (which holds for $i$, since it follows from Eq. \eqref{eqchainenv} by taking $K = H_{i-1}$), 
 we obtain $ \chi(K) < 0$; this   completes the proof of \eqref{eqchi2} and that of Eq. \eqref{eqchaine2} for any $K$.

Now let $K$ be a \cas of $G$ such that $\dim K > \dimhi$ and  $\chi(K) = 0$.  Let us prove that $H_i \subset K$ and $K \subset H_{i+1}$; this will conclude the proofs of the equality cases in Eqns.  \eqref{eqchaine1} and  \eqref{eqchaine2}, and will also prove that $H_i  \subset H_{i+1}$ since one may take $K=H_{i+1}$. 

With this aim in view, we notice, using \eqref{eqchi1}, 
that  $ \dim K \leq \dimhipu$  and that for any $K'$ we have
\begin{equation} \label{eqsommeinter3}
\chi(K \cap K')+\chi(K+K') \geq  \chi(K) + \chi(K')
\end{equation}
using Eqns. \eqref{eqsommeinter1} and  \eqref{eqsommeinter2}. Recall that $\chi(H_i) = \chi(H_{i+1})= \chi(K)=0$ thanks to Eq. \eqref{eq8bis} and our assumption on $K$,  and that $\chi(K'') \leq 0$ for any $K''$ using \eqref{eqchi1} and \eqref{eqchi2}. With $K' = H_i$ we obtain in this way $\chi(K\cap H_i) = \chi(K+ H_i) = 0$, so that   \eqref{eqchi2} yields $K\cap H_i = H_i$  and finally $H_i \subset K$. In a similar way, with $K'=H_{i+1}$ we get $\chi(K\cap H_{i+1}) = \chi(K+ H_{i+1}) = 0$, so that   \eqref{eqchi1}  yields $\dim(K+H_{i+1})=\dimhipu$  and thus $K \subset H_{i+1}$. This concludes the proof that $H_i \subset K \subset H_{i+1}$, and that of Proposition \ref{propchainez}   and Remark  \ref{remchainei}.

\subsection{Independence and finiteness results} \label{subsec320}

It is clear from the construction that the subgroups $H_0$, \ldots, $H_r$ and the integer $r$ depend on $\Sul$. However there is a transformation under which they are invariant:

\begin{Lemme}  \label{lemchaineiv} The subgroups $H_0$, \ldots, $H_r$ and the integer $r$ remain the same if $\Sul$ are replaced with $S_1^\alpha$, \ldots,  $S_l^\alpha$ for some $\alpha>0$.
\end{Lemme}

An important consequence of this lemma is that if $S_1=\ldots=S_l = S > 1$ as in the introduction, then  $H_0$, \ldots, $H_r$ do not depend on $S$ (see \S \ref{subsec33}). 

\bigskip

\Dem of Lemma \ref{lemchaineiv}: 
Upon replacing $\Sul$  with  $S_1^\alpha$, \ldots,  $S_l^\alpha$, the function $\fis$ is multiplied by $\alpha>0$ so that Eq.  \eqref{eqchaine1} is still valid:  since the   chain of subgroups constructed above is unique (see Proposition \ref{propchainez}), it remains the same.

\bigskip

Throughout the proof of Theorems \ref{th11} and \ref{th12}, many constants will appear that depend on the subgroups $\chaine$. The following lemma shows that such a constant can be made independent from these subgroups, by increasing it if necessary.

\begin{Lemme}     \label{lemchainev}  There exists a finite set $\calE$, which depends on $\gamul$ but not on  $\Sul$,  such that   all subgroups $H_0$, \ldots, $H_r$ belong to $\calE$.
\end{Lemme}

The idea behind this lemma is simply that the construction of $H_i$ involves only  $\dim H_i$ and the ranks of $\Gamma_j\cap H_i$, which take only finitely many values (see  \S 3.1.1 of \cite{SFinterp} for an analogous  situation).  Moreover there is no \cas $H' \neq H_i$ such that $\dim H' = \dimhi $ and $\rk ( \Gamma_j\cap H' ) = \rk ( \Gamma_j\cap H_i )$ for any $j$, so that $H_i$ can take only finitely many values. In the notation of  Remark \ref{rempolygone}, there is no $H' \neq H_i$  such that $M_{H'} = M_{H_i}$, even
though in general there may exist \cass $H' \neq H''$ such that $M_{H'} = M_{H''}$: the  \cass $H_i$ are uniquely determined by the vertices of the polygon $M_{H_0}M_{H_1}\ldots M_{H_r}N$.
Let us make these ideas more precise now.

\bigskip

\Dem of Lemma \ref{lemchainev}:  Let us denote by $\calS$ the set of all \cass of $G$, and for $K \in \calS$ we let $\psi(K) = (\dim K, \rk(\Gamma_1\cap K), \ldots,  \rk(\Gamma_l\cap K))\in\Z^{l+1}$. 
Let $\calE $ denote the set of all $H \in \calS$ such that $\psi^{-1}(\psi(H))= \{H\}$. Then $\calE$ is a finite set, because  $\psi(\calS)$ clearly is and $\psi_{|\calE} : \calE \to \psi(\calS)$ is an injective map.  Now  for any subgroup  $H_{i+1}$  in a chain corresponding to some $\Sul$, we have $\psi^{-1}(\psi(H_{i+1}))= \{H_{i+1}\}$ because equality holds in Eq.  \eqref{eqchaine1}    for any $K \in \psi^{-1}(\psi(H_{i+1}))$; therefore $H_{i+1}\in\calE$.  Since $\{0\} \in\calE$, this concludes the proof of Lemma \ref{lemchainev}.

\subsection{Applications to the Distribution of $\Gamma$} \label{subsec32}

The chain of algebraic subgroups $\chaine$ constructed in \S \ref{subsec31} is useful to study the distribution of $\GS$ with respect to algebraic subgroups of $G$. Several results of this kind have been stated in the introduction when $S_1=\ldots=S_l$, and will be proved in \S \ref{subsec33}. In the general case where $\Sul$ are not assumed to be equal, the same results hold except that they have to be stated differently: exponents like $\mu(\Gamma,G)$ and $\muet(\Gamma,G)$ are no  longer available. We shall neither state nor prove the corresponding generalizations of all results stated in the introduction, but only the ones that will be used in the proof of Theorems \ref{th11} and \ref{th12}. 

To begin with, let us generalize the fact that $\Gamma(S) \cap H_{i+1} \bmod H_i$ is  well-distributed in  $H_{i+1}/H_i$. Recall that $\sss_i$ has been defined at the beginning of \S \ref{sec1}.

\begin{Lemme} \label{lem377} There exists a positive constant $\ctss$, which depends only on $G$, $\gamul$ but not on $\Sul$, with the following property: for any $i \in \zerormu$ and any \cas $H$ such that $H_i \subsetneq H$, we have
$$\Card\Big(\Gamma(2n\Ssoul)\cap H \bmod H_i\Big) < \ctss \sss_i ^{\dim ( H/H_i)}$$
where $n = \dim G$.
\end{Lemme}

This lemma asserts that in applying Proposition \ref{propinterp} to $\Gamma(2n\Ssoul)\cap H_{i+1} \bmod H_i$ in the algebraic group $H_{i+1}/H_i$, it is enough to check 
assumption \eqref{eqhyppropinterp} with $H = H_{i+1}/H_i$ (with a smaller value of $\cinterp$, though), so that this proposition applies as soon as $D > \cnvu \sss_i T$ for some constant $\cnvu$. Indeed $\sss_i$ is  given by 
\begin{equation} \label{eqgensi}
\sss_i = \exp\Big[\frac{\fis(H_{i+1})-\fis(H_i)}{\dimhipu-\dimhi}\Big],
\end{equation}
and Lemma \ref{lemab} shows that $\sss_i^{\dimhipu-\dimhi}$ is equal, up to a multiplicative constant, to the cardinality of $(\GS \cap H_{i+1} )  \bmod H_i$: the conclusion of Lemma \ref{lem377} is an equality for $H= H_{i+1}$, except for the value of the constant $ \ctss$.

We prove Lemma \ref{lem377} for $\Gamma(2n\Ssoul)$ because it will be applied in this way in the proof of Theorem  \ref{th12}. The value $2n$ could be replaced with any  other constant $\cnvt$ and then
 $\ctss$ would depend on $\cnvt$.  Notice also  that the chain of algebraic subgroups associated (as in  \S \ref{subsec31}) with the parameters $2nS_1,\ldots,2nS_l$ might be distinct from
the chain $\chaine$  associated with $S_1,\ldots,S_l$ (which appears in  Lemma \ref{lem377}). 

\bigskip

\Dem of Lemma \ref{lem377}: Lemma \ref{lemab} applied to $\Gamma(2n\Ssoul)$, $H'=H$ and $H''=H_i$ yields
$$\Card\Big(\Gamma(2n\Ssoul)\cap H \bmod H_i\Big) < \cabbpr \exp\Big(\fis(H)-\fis(H_i)\Big)$$
where $ \cabbpr $ depends only on $\gamul$ and $n$, using Remark \ref{remabd} and Lemma \ref{lemchainev}. Since
$$\sss_i^{\dim(H/H_i)} = \exp \Big[ \frac{\dim H - \dimhi}{\dimhipu - \dimhi}\Big(\fis(H_{i+1})-\fis(H_i)\Big)\Big],$$
Lemma \ref{lem377} follows using Eq. \eqref{eqchaine1} of Proposition \ref{propchainez}.

\bigskip

The next lemma corresponds, when $S_1=\ldots=S_l$,  to the result $\mu(\Gamma\cap H_i , H_i) = \frac{\rk(\Gamma\cap H_i)-\rk(\Gamma\cap H_{i-1})}{\dimhi - \dimhimu}$.

\begin{Lemme} \label{lem373} For any $\eps  > 0$ there exists a positive constant $\ctst$, which depends only on $\eps$, $G$, $\gamul$ but not on $\Sul$, with the following property: for any $i \in \unr$ and any \cas $H$ such that $H\subsetneq H_i$, we have
$$\Card\Big(\Gamma\Big(\frac{\eps}{\dim H_i}\Ssoul\Big)\cap H_i \bmod H \Big) > \ctst \sss_{i-1} ^{\dim ( H_i/H)}.$$
\end{Lemme}

This lemma asserts that in applying Proposition \ref{proppphz} to $\Omega = \Gamma(\frac{\eps}{\dim H_i}\Ssoul)\cap  H_i$ in the algebraic group $H_i$, it is enough to check 
assumption \eqref{eqhypproppphz} with $H = H_{i-1}$ (with a larger value of $\cpphz$, though), so that this proposition applies as soon as $D < \cnvd \sss_{i-1} T$ for some constant $\cnvd$. Indeed 
$\sss_{i-1}^{\dim(H_i/H_{i-1}) }$ is equal, up to a multiplicative constant, to the cardinality of $( \Gamma(\frac{\eps}{\dim H_i}\Ssoul) \cap H_{i } )  \bmod H_{i-1}$: the conclusion of Lemma \ref{lem373} is an equality for $H=H_{i-1}$, up to the value of $\ctst$.

As for Lemma \ref{lem377} above,
we prove Lemma \ref{lem373} for $\Gamma(\frac{\eps}{\dim H_i}\Ssoul)$ because it will be applied in this way in the proof of Theorem  \ref{th11}. The value $\frac{\eps}{\dim H_i}$ could be replaced with another constant. 

\bigskip

\Dem of Lemma \ref{lem373}: Applying Lemma \ref{lemab}   to $\Gamma\Big(\frac{\eps}{\dim H_i}\Ssoul\Big) $ with $H'=H_i$ and $H''=H$ yields, using Remark \ref{remabd} and Lemma \ref{lemchainev}:
$$\Card\Big(\Gamma\Big(\frac{\eps}{\dim H_i}\Ssoul\Big)\cap H_i \bmod H \Big) >  \cabapr \exp\Big(\fis(H_i)-\fis(H)\Big)$$
where $ \cabapr $ depends only on $\gamul$ and $\eps$. Since
$$\sss_{i-1}^{\dim(H_i/H )} = \exp \Big[ \frac{\dim H_i - \dim H}{\dimhi  - \dimhimu}\Big(\fis(H_i)-\fis(H_{i-1})\Big)\Big],$$
Eq. \eqref{eqchaine3} in Proposition \ref{propchaineii}  enables one to conclude the proof of Lemma \ref{lem377}.

\subsection{The case $S_1=\ldots=S_l$}\label{subsec33}

In this subsection we assume that  $S_1=\ldots=S_l  = S > 1 $ as in the introduction. We shall deduce the results announced there  from those proved previously in the general setting.

To begin with, we have $\fis(K) = \rk(\Gamma\cap K) \log S$ for any \cas $K$ of $G$. In particular $\fis(K)$ does not depend on $\gamul$, but only on $\Gamma$, $S$ and $K$. It is easily seen that the factor $\log S$ cancels out in all inequalities like Eq. \eqref{eqchaine1}, so that the chain of algebraic subgroups $\chaine$ depends only on $\Gamma$ but not on $\gamul$ or on $S$ (the latter point  is also a consequence of Lemma \ref{lemchaineiv}). 

On the other hand  Eq. \eqref{eqchaine1} is equivalent to 
$$ \frac{\rk(\Gamma \cap K) - \rk(\Gamma \cap H_i)}{\dim K - \dim H_i} \leq \frac{\rk(\Gamma \cap H_{i+1}) -\rk(\Gamma \cap H_i)}{\dim H_{i+1} - \dim H_i}.$$
Since $H_0 =\{0\}$, Proposition \ref{propchainez} asserts that $\frac{\rk(\Gamma \cap K)}{\dim K} \leq \frac{\rk(\Gamma \cap H_1)}{\dim H_1}$, 
and if equality holds then $ K \subset H_{1}$. This implies $ \muet(\Gamma,G) = \frac{ \rk(\Gamma\cap H_1)}{  \dim H_1}$.
In the same way, $H_r = G$ so that Proposition \ref{propchaineii} yields $\frac{\rk(\Gamma) - \rk(\Gamma \cap K)}{\dim G  - \dim K} \geq \frac{\rk(\Gamma) - \rk(\Gamma \cap H_{r-1})}{\dim G  - \dim H_{r-1}}$, 
and if equality holds then $ H_{r-1}   \subset K$. In particular we have  $ \mu (\Gamma,G) = \frac{\rk(\Gamma) - \rk(\Gamma\cap H_{r-1})}{ \dim G  - \dim H_{r-1}}$.

Letting $\mu_i  =  \frac{\rk(\Gamma\cap H_{i+1}) - \rk(\Gamma\cap H_{i})}{\dim H_{i+1} - \dim H_i}$ as in the introduction, with $i\in\zerormu$, Proposition \ref{propchaineiii} asserts that $\mu_{r-1} < \ldots < \mu_1< \mu_0$; in the notation of Remark \ref{rempolygone} this is clear because $\mu_i \log S $ is the slope of the line $(M_{H_i}M_{H_{i+1}})$.  Lemma \ref{lemab} proves that $\Card (( \GS \cap H_{i+1})\bmod H_i)$ is equal, up to a multiplicative constant depending only on $\gamul$, to $S^{{\rk(\Gamma\cap H_{i+1}) - \rk(\Gamma\cap H_{i})}} = S^{\mu_i ( \dim H_{i+1} - \dim H_i) }$; the quantity denoted by $\sss_i$ in this paper equals $S^{\mu_i}$ in this case.

Let us fix $i$, $j$ such that  $0\leq i < j \leq r$. Then the chain of algebraic subgroups associated with $\Gamma \cap H_j \bmod H_i$ in the algebraic group $H_j/H_i$  is  $\{0\} = \frac{H_i}{H_i}  \subsetneq \frac{H_{i+1}}{H_i} \subsetneq \ldots  \subsetneq \frac{H_{j-1}}{H_i}  \subsetneq  \frac{H_j}{H_i}$. This proves the equalities
$$\mu(\Gamma \cap H_j \bmod H_i , H_j/H_i) = \frac{\rk(\Gamma\cap H_j) - \rk(\Gamma\cap H_{j-1})}{\dim H_j - \dim H_{j-1}}$$
and
$$\muet(\Gamma \cap H_j \bmod H_i , H_j/H_i) = \frac{\rk\Big(\frac{ \Gamma\cap H_{i+1}}{ \Gamma\cap H_i}\Big)}{\dim ( H_{i+1}/H_i)}.$$

\section{Proof of the Main Result} \label{sec4}

In this section we prove Theorems \ref{th11} and \ref{th12}. The strategy is to apply the special cases where $D$ is very large or very small, proved in \S \ref{subsec21}, to sub-quotients of $G$. We deduce Theorem \ref{th12} from a multiplicity estimate in the algebraic group $H_i$. The proof of Theorem \ref{th11} is more complicated: it involves   interpolation estimates in $H_{i+1}/H_i$  and in $G/H_{i+1}$.

\bigskip

The addition law and the translations on $G$  play a key role in the proof and we begin by
establishing our notation to represent these.
Let $a$ and $b$ be integers such that 
there exists a complete system of addition laws on $G$ of bi-degree
$(a,b)$. This means that the addition law on $G$ (embedded in $\P^N$)
is represented, on every element of a suitable open cover, by a family
of bi-homogeneous polynomials of bi-degree $(a,b)$.  We may assume (see \cite{MWun}, p. 493) that some open set $U$ in this cover contains  $\Gamma$ and the point $\gamma$ introduced below at the beginning of the proof of  Theorem \ref{th11}. Of course $U$ depends on $\gamma$, but this is not important in the proof. There exists a family $E_0(X,Y)$, \ldots, $E_N(X,Y)$ of bi-homogeneous polynomials of bi-degree $(a,b)$ which represents the addition law on $U \times U$; we let $X = (X_0,\ldots,X_N)$, $Y = (Y_0,\ldots,Y_N)$ and  $E = (E_0,\ldots,E_N)$.

For any $y\in U$, after choosing a system $(y_0,\ldots,y_N) \in \C^{N+1}$ of projective coordinates of  $y$ in $\P^N$ we may consider for any $P \in \cdegD$ the polynomial
$$t_yP(X) = P(E(X,y))\in\cdeg_{aD}.$$
The linear map $t_y : \cdegD \to \cdeg_{aD}$ represents the translation by $y$.  Moreover
 if $P$ vanishes to order at least $T$ (resp. does not vanish) at a given point $z$ and if $z-y\in U$ then $t_y P$ vanishes to order at least $T$ (resp. does not vanish) at the point $z-y$.
 Of course the map $t_y$ depends on $D$, $E$ and on the choice of  $(y_0,\ldots,y_N)$, but we omit this dependence in the notation $t_y$.

\bigskip

In the proof we shall use repeatedly the following fact: since $H_i$ may take only finitely many values (see Lemma \ref{lemchainev}), a constant that depends on $H_i$ can actually be chosen in terms of $G$, $\gamul$. Given a constant $N$ (depending on $\gamul$), we may also assume that $D$ is a multiple of $N$. Indeed   $\base$ is a non-increasing function of $D$ when the subset $\GS$ and the order of vanishing $T$ are held constant, so it is enough to prove Theorem \ref{th11} for a slightly smaller value of $D$ (resp. to prove Theorem \ref{th12} for a slightly larger  value of $D$).

\bigskip

We first prove   Theorem \ref{th12}. Let $\gamma \in \Gamma((1-\eps)\Ssoul)$ and assume that $P \in \cdegD$  vanishes to order at least $T$ at any point of $\GS$. Consider $Q = t_\gamma P \in \cdeg_{aD}$: then $Q$  vanishes to order at least $T$ at any point of $ \Gamma(  \eps \Ssoul)$. We let $\Omega_1 =  \Gamma( \frac{ \eps }{\dimhi} \Ssoul) \cap H_i$ and denote by $Q_1 \in\cdehi_{aD}$ the restriction of $Q$ to $H_{i}$. Then $Q_1$  vanishes to order at least $T$ at any point of  $\Omega_1[\dimhi]$. Moreover, for any \cas $H \subsetneq H_i$ Lemma \ref{lem373} yields
$$\Card (\Omega_1 \bmod H ) > \ctst \sss_{i-1}^{\dim(H_i/H)}$$
where $\ctst$ depends only on $\gamul$. Since $D < \cthud^{-1} \sss_{i-1}T$ this implies (provided that $\cthud$ is large enough) that
$$\Card (\Omega_1 \bmod H ) T ^{\dim(H_i/H)} > \cpphz D^{\dim(H_i/H)}$$
where $\cpphz$ is the constant in Proposition \ref{proppphz} applied in the algebraic group   $  H_i$. This Proposition yields $Q_1=0 \in \cdehi_{aD}$ so that $P $ vanishes identically on $\gamma + (H_i \cap U)$.  Now the zero element of $G$ belongs to $  H_i \cap U$ (because we have assumed $\Gamma \subset U$), so that $H_i \cap U$ is non-empty. Since $U$ is an open subset of $G$, we obtain that $H_i \cap U$ is Zariski dense in $H_i$. This density does not change by translation, so that   $P $ vanishes identically on  $\gamma + H_i$ because it vanishes
 on $\gamma + (H_i \cap U)$. This   concludes the proof of  Theorem \ref{th12}.  
 
\bigskip

We now prove   Theorem \ref{th11}. We argue by decreasing induction on $i$. Letting $\sss_r = 1$ this result is meaningful for $i=r$, and trivially true since $H_r=G$. From now on, we let $i\in\zerormu$ and assume that  Theorem \ref{th11}  holds for $i+1$.

Assume there exists $\gamma\in\base$ with $\gamma\not\in \GS+H_i$. Since    Theorem \ref{th11}  holds for $i+1$ and $\sss_{i+1}\leq\sss_i$, we have $\gamma \in \GS+H_{i+1}$. Let $\beta\in\GS$ and $h\in H_{i+1}$ be such that $\gamma = \beta+h$. Consider $\Omega_2 = (-\beta+\GS)\cap H_{i+1}$, and notice that $h \not\in\Omega_2+H_i$ since $\gamma\not\in \GS+H_i$. 

Now $H_{i+1}/H_i$ is a  commutative algebraic group, so we can choose (arbitrarily) a projective embedding $H_{i+1}/H_i \hookrightarrow \P^{M_i}$. With respect to this embedding (and that of $H_{i+1}$ in $\P^N$), the projection $H_{i+1} \to H_{i+1}/H_i $ is given, on an open subset of $H_{i+1}$ which contains $\Omega_2\cup\{h\}$,  by homogeneous polynomials $R_{i,0}$, \ldots, $R_{i,M_i}$ of the same degree, say $a_i$. It is possible to ensure that $a_i$ depends only on the embeddings of $H_{i+1}/H_i$ and $H_{i+1}$, and not on $\Omega_2$ or $h$. We put $R_i = (R_{i,0}, \ldots,  R_{i,M_i})$. 

Let $\hbar = h \bmod H_i$ and $\Omegabar_2 = \Omega_2 \bmod H_i$, so that $\hbar\not\in\Omegabar_2$. Then $(\Omegabar_2\setminus\{\hbar\})\{\dim(H_{i+1}/H_i)\}$ is a subset of $\Gamma(2n\Ssoul)\cap H_{i+1} \bmod H_i$ since $\dim(H_{i+1}/H_i)\leq \dim G = n$, so that Lemma \ref{lem377} yields
$$\Card \Big( (\Omegabar_2\setminus\{\hbar\})\{ \dim(H_{i+1}/H_i)\} \cap \Hbar\Big) < \ctss \sss_i^{\dim\Hbar}$$
for any \cas $H$ of $G$ such that $H_i \subsetneq H \subset H_{i+1}$, where $\Hbar = H/H_i$. Since $D > \cthuu \sss_i T$, Proposition \ref{propinterp} applies in  the algebraic group $  H_{i+1}/H_i$ if $\cthuu$ is sufficiently large: it provides $P_1 \in \cdequo_{D/2aa_i}$ which  vanishes to order at least $T$ at any point of  $\Omegabar_2$ and does not vanish at $\hbar$ (because $\hbar\not\in \Omegabar_2$). 
Then $P_1 \circ R_i \in \cdehipu_{D/2a}$ vanishes to order at least $T$ at any point of  $\Omega _2 = (-\beta+\GS)\cap H_{i+1}$, and does not vanish at  the point $ h$. Choose $P_2 \in \cdeg_{D/2a}$ such that $ P_1 \circ R_i $ is the restriction of $P_2$ to $H_{i+1}$
so that the same vanishing and non-vanishing properties hold for $P_2$. Therefore $P_3 = t_{-\beta}P_2\in\cdeg_{D/2}$ vanishes to order at least $T$ at any point of  $ \GS \cap (  \beta+ H_{i+1}) $, and does not vanish at  the point $\beta+ h = \gamma$ (because $U$ contains $\Gamma$ and $\gamma$).

On the other hand, we can choose an embedding of the  commutative algebraic group $G/H_{i+1}$ in   projective space $\P^{M'_i}$.
On an open subset of $G$  which contains $\Gamma \cup\{ \gamma\}$ the   projection $G\to G / H_{i+1} $ is given as above by a family $R'_i = (R'_{i,0}, \ldots,  R'_{i,M'_i})$ of homogeneous polynomials  of the same degree  $a'_i $ and this degree can be made independent from $\gamma$. 

Now let $\Omegabar_3 = \overline{\GS} \setminus\{\betabar\}$, where $ \overline{\GS}$ and $\betabar$ are the images of $\GS$ and $\beta$ in $G/H_{i+1}$. Then $\Omegabar_3\{ \dim(G/H_{i+1})\} \subset \Gamma(2n\Ssoul)\bmod H_{i+1}$ since $\dim(G/H_{i+1})\leq n$, so that Lemma \ref{lem377} (with $i+1$ instead of $i$) yields
$$\Card \Big( \Omegabar_3 \{ \dim(G/H_{i+1} ) \} \cap \Hbar \Big) < \ctss \sss_{i+1}^{\dim\Hbar}$$
for any \cas $H$ such that $H_{i+1} \subsetneq H \subset G$, with $\Hbar = H/H_{i+1}$. Now we have $D > \cthuu \sss_i T \geq \cthuu \sss_{i+1}T$; if $\cthuu$ is sufficiently large then Proposition \ref{propinterp} (applied in $G/H_{i+1}$) provides $Q_1 \in\cdegsh_{D/2a'_i}$ which vanishes  to order at least $T$ at any point of  $\Omegabar_3 = \overline{\GS} \setminus\{\betabar\}$, and does not vanish at  the point $\betabar$. Then $Q_2 = Q_1\circ R'_i \in\cdeg_{D/2}$ vanishes  to order at least $T$ at any point of  $\GS \setminus ( \beta+H_{i+1})$, and does not vanish at $\gamma = \beta+h$.

We consider now $P = P_3Q_2 \in\cdegD$. We have $P(\gamma)\neq 0$ because $P_3(\gamma)\neq 0$ and $Q_2(\gamma)\neq 0$.  For $\gamma'\in\GS$,  if $\gamma'\in\beta+H_{i+1}$ then $P_3$ vanishes  to order at least $T$ at $\gamma'$; otherwise $Q_2$ does. Therefore $P$ vanishes  to order at least $T$ at any point of $\GS$. Since $P(\gamma)\neq 0$ and $\gamma\in\base$, this is a contradiction.

\section{Possible generalizations} \label{secgen}

It would be interesting to generalize Theorems \ref{th11} and \ref{th12}  in at least two directions.

The first one would be to replace $\GS$ with a fixed finite set $\Omega$.  A first step   would be   to find   constants $\cthuu$ and $\cthud$ in Theorems \ref{th11} and \ref{th12}  which do  not depend on $\gamul$. Interpolation and multiplicity estimates are known in this setting (see \S \ref{subsec21}), but the proof  leaves no hope to obtain this result unless new ideas are used. For instance, the chain of subgroups $\chaine$ constructed in \S \ref{sec3} does not depend on the torsion part of $\Gamma$: it is trivial as soon as $\Gamma$ has rank 0. For an analogous reason, in Masser's interpolation estimate \cite{Masser} (and in the first author's generalization \cite{SFinterp}), the constant depends also on $\gamul$. The case of an arbitrary  finite set $\Omega$ was dealt with in \cite{FNinterp} using a more geometric approach in terms of Seshadri constants.

The second way to generalize Theorems \ref{th11} and \ref{th12} would be to consider vanishing along analytic subgroups of $G$. The only problem is that the interpolation estimate of \cite{FNinterp} (stated above as Proposition  \ref{propinterpFN}) is not known in this setting. The only available interpolation estimate is the one proved by the first author  \cite{SFinterp}, but it is not sufficiently precise to deduce the corresponding generalization of Proposition \ref{propinterp} (even with $\Omega = \GS$). 

Let us introduce this setting more precisely, and mention   how the other 
tools used in this paper generalize to it.
Let  $T_0(G)$  denote  the tangent space to $G$ at  0,  seen as   the space of translation-invariant vector fields on $G$.
Let $W$ be a subspace of  $T_0(G)$, of dimension $d \geq 0$, and $\drondsoul = (\drond_1,\ldots,\drond_d)$ be a basis of  $W$. 
For a family $\Tsoul = (T_1,\ldots,T_d)$ of $d$ positive real numbers, 
we let $\Ndt$ be the set of all 
$\sigma = (\sigma_1,\ldots,\sigma_d) $ such that $0 \leq \sigma_j 
< T_j$ for any $j \in \und$. 
  
We denote by $\opW$ the set of all polynomials in $\drond_1,\ldots,\drond_d$,
i.e. the space of differential operators along $W$, and   by 
$\opdrondT$ the subspace of $\opW$ spanned by the monomials
$\drondsoul^\sigma = \drond_1 ^{\sigma_1} \ldots \drond_d ^{\sigma_d}$ for
$\sigma \in \Ndt$. We assume (without loss of generality: see \cite{MWun}, p. 492) that 
 $
\Gamma \inclus \{X_0 \neq 0 \} \inclus \P^N 
 $, and we say that
 a polynomial $P \in \cdegD$ vanishes up to order
$\Tsoul$ along $W$ at a point $\gamma \in \Gamma$ if 
$\drond ^\sigma (P/\xzd) (\gamma) =
0$ for any $\sigma \in \Ndt$. If $T_1=\ldots=T_d=T$ and $W =T_0(G)$, this means that $P$ vanishes to order at least $T$ at $\gamma$.

With this notation, one would replace everywhere ``vanishing to order at least $T$'' with ``vanishing up to order
$\Tsoul$ along $W$'', and $T^{\dim H}$ with  $\dim (\opWTH \inter \opdrondT)$. The corresponding multiplicity  estimate (Proposition \ref{proppphz}) has been proved by Philippon  \cite{PphRocky}.
For any  $j\in\{0,\ldots,d\}$ let $W_j = \Span(\drond_1,\ldots,\drond_j)$.  Then $\dim(W_j\cap T_0(H))$ plays the role of $\rk(\Gamma_j\cap H)$. Given $\Sulgeq \geq 1$ and $T_1\geq\ldots\geq T_d\geq 1$, let us define
$$\fist(K) = \sum_{j=1}^l \rk\Big(\frac{\Gamma_j\cap K}{\Gamma_{j-1}\cap K}\Big) \log S_j +  \sum_{j=1}^d \dim\Big(\frac{W_j\cap T_0(K)}{W_{j-1}\cap T_0(K)}\Big) \log T_j$$
for any \cas $K$ of $G$. Then in \S \ref{subsec31} it is possible to replace $\fis$ with $\fist$.  
The chain of subgroups constructed in this way depends on $\gamul$, $\drond_1,\ldots,\drond_d$, $\Sul$, $T_1,\ldots,T_d$. In Lemmas \ref{lem377} and \ref{lem373}, the left hand 
side of the inequalities has to be multiplied by $\dim\Big(\frac{  \opWTHpr \inter \opdrondT}{  \opWTHsec \inter \opdrondT}\Big)$ with $\{H',H''\} = \{H,H_i\}$. Moreover in the definition of $\sss_i$, $\fis$ should also be replaced with $\fist$ (see Eq. \eqref{eqgensi}). It seems reasonable to conjecture that Theorems \ref{th11} and \ref{th12}  hold  in this setting, with $\sss_i T$ replaced with this new value of $\sss_i$.

 \newcommand{\url}{\texttt}

\providecommand{\bysame}{\leavevmode ---\ }
\providecommand{\og}{``}
\providecommand{\fg}{''}
\providecommand{\smfandname}{\&}
\providecommand{\smfedsname}{eds.}
\providecommand{\smfedname}{ed.}
\providecommand{\smfmastersthesisname}{M\'emoire}
\providecommand{\smfphdthesisname}{Th\`ese}

\end{document}